\documentclass{ws-ijbc}
\usepackage[most]{tcolorbox}
\usepackage{subfloat}
\usepackage{graphicx}
\usepackage{subfigure}
\usepackage{epstopdf}
\usepackage{amssymb} 
\usepackage{longtable}
\usepackage{algorithm}
\usepackage{algorithmic}
\usepackage{graphicx}
\usepackage{pgf, tikz}
\usetikzlibrary{arrows, automata}
\usepackage{amsmath}
\usepackage{lipsum}
\usepackage{float}
\usepackage{multicol}
\usepackage{longtable}

\newtheorem{thm}{Theorem}
\newtheorem{conj}{Conjecture}

\newtheorem{defn}{Definition}

\newtheorem{pro}{Proposition}
\usepackage{multirow}
\usepackage{tkz-graph}

\begin{document}

%\catchline{}{}{}{}{} % Publisher's Area please ignore

%\markboth{Author's Name}{Paper Title}

\title{Two Dimensional Discrete Dynamics of Integral Value Transformations}

\author{Jayanta Kumar Das}

\address{Applied Statistics Unit, Indian Statistical Institute\\ Kolkata-$700108$, India\\
dasjayantakumar89@gmail.com}

\author{Sudhakar Sahoo}
\address{Institute of Mathematics and Applications\\ Bhubaneswar-$751029$,  India.\\
sudhakar.sahoo@gmail.com}

\author{Sk. Sarif Hassan\footnote{Corresponding Author}}
\address{Department of Mathematics, Pingla Thana Mahavidyalaya\\
Paschim Medinipur-$721140$,  India.\\
sarimif@gmail.com}

\author{Pabitra Pal Choudhury}
\address{Applied Statistics Unit, Indian Statistical Institute\\ Kolkata-$700108$, India\\
pabitrapalchoudhury@gmail.com}

\maketitle

\begin{history}
\received{(to be inserted by publisher)}
\end{history}

\begin{abstract}
A notion of dimension preservative map, \textit{Integral Value Transformations} (IVTs) is defined over $\mathbb{N}^k$ using the set of $p$-adic functions. Thereafter, two dimensional \textit{Integral Value Transformations} (IVTs) is systematically analyzed over $\mathbb{N} \times \mathbb{N}$ using pair of two variable Boolean functions. The dynamics of IVTs over $\mathbb{N} \times \mathbb{N}=\mathbb{N}^2$ is studied from algebraic perspective. It is seen that the dynamics of the IVTs solely depends on the dynamics (state transition diagram) of the pair of two variable Boolean functions. A set of sixteen \textit{Collatz-like} IVTs are identified in two dimensions. Also, the dynamical system of IVTs having attractor with one, two, three and four cycles are studied. Additionally, some quantitative information of \textit{Integral Value Transformations} (IVTs) in different bases and dimensions are also discussed. 
\end{abstract}


%\begin{multicols}{2}
\section{Introduction}
\noindent
Classically, discrete dynamics refers to the study of the iteration of self-maps of the complex plane or real line \cite{1,2,3,22}. There are plenty of works carried out over last couple of decades in discrete dynamical systems as well as in arithmetic dynamics \cite{4,5,6,7,8,21,23,24,25,26,27}. What is not very well studied yet is the discrete dynamics over spaces like natural numbers. In $2010$, the notion of Integral Value Transformations (IVTs) which is a new paradigm in arithmetic dynamics \cite{9,10,11} has been introduced. Later in $2014$, a special class of one dimensional Integral Value Transformations (called \textit{Collatz-like}) is studied in \cite{12}. A couple of other works on IVTs and their applications in arithmetic logical circuit design and in Biology are also studied \cite{13,14,15,16,17,18,19}. Applications of such transformations are noteworthy in the area of Cryptology and other field of discrete mathematics \cite{29,30,31,32,33}.

In this present article, an attempt has been initiated to understand the discrete dynamics of Collatz-like IVTs including other two dimensional Integral Value Transformations. It is seen that the dynamics of the IVTs are depending on the state transition diagram of two dimensional IVTs of the pair of two variable Boolean functions \cite{20}. There are four types of dynamics seen with regards to attractors with different cycle length(s) $1, 2, 3$ and $4$. All these dynamical systems end with attractors (global and/or local) with a maximum of four length cycles. In addition, from the algebraic aspects the dynamical system of IVTs are investigated. Discrete dynamics is itself very appealing from its application viewpoint. The present study of the dynamical systems over two dimensions can be as pathways for network designing converging to the desired attractor(s).       

The rest of paper is organized as follows. In Section $2$, the definitions of IVT in higher dimensions and in particular in two dimensions over the set of natural numbers with zero are given. The various algebraic properties of the dynamical system of IVTs are described in this section. The Section $3$ deals with the state transition diagram for the IVTs of various algebraic classes. Some important theorems including all possible type of attractors are presented in this section. The Section $4$ deals with the dynamics of IVTs using their profound trajectories. The Section $4$ also deals with the classifications of $256$ possible IVTs in two dimensions based on their dynamics over $\mathbb{N}\times \mathbb{N}$. The quantitative measure for the various classes in different bases and dimensions are also incorporated in this section. This Section 4 also includes the dynamics of some special class of IVTs. The Section $5$ concludes the key findings of current article and future research endeavors.

\section{Definitions \& Basic Algebraic Structures in $\mathbb{B}_2$ and $\mathbb{T}$}
\subsection{Definitions and basic algebraic structures}

\begin{defn} 
\emph{Integral Value Transformation in the base $p$ and dimension $k$} $(IVT^{p, k}_{i_{1},i_{2},\cdots ,i_{k}})$ over $\mathbb{N}^k$ ($\mathbb{N}$ is set of all natural numbers including zero) is defined as 
\begin{equation}
\begin{split}
IVT^{p,k}_{i_{1},i_{2},\cdots ,i_{k}}(m_{1},m_{2},\cdots, m_{k})=(IVT^{p,k}_{i_{1}} (m_{1},m_{2},\cdots, m_{k}), IVT^{p,k}_{i_{2}} (m_{1},m_{2},\cdots, m_{k}),\cdots,\\ IVT^{p,k}_{i_{k}} (m_{1},m_{2},\cdots, m_{k})) \nonumber
\end{split} 
\end{equation}
	
\noindent
Here $(m_{1},m_{2},\cdots, m_{k})$ is an element of $\mathbb{N}^k$. The function $IVT^{p,k}_i$ is defined using $f_i$ in \cite{11}. $f_i$ is a $p$-adic function. 
\end{defn}

\begin{defn} 
\emph{Integral Value Transformation in two dimension (k=2) with base (p=2):} $(IVT^{2,2}_{i,j})$ over $\mathbb{N} \times \mathbb{N}$ is defined as 
	$$IVT^{2,2}_{i,j}(m,n)=(IVT^{2,2}_i (m,n), IVT^{2,2}_j (m,n))$$
	
\noindent
Here $(m,n)$ is an element of $\mathbb{N} \times \mathbb{N}$ and the function $f_i$ or $f_{j}$ is the $2$-variable Boolean functions and pair-wise in combining they are represented as $f_{i,j}$ which is formally defined in the Definition 3.

\end{defn}

\begin{defn}
The map $f_{i,j}$ is defined on the set $\mathbb{B}_2=\{00,01,10,11\}$ as 
$$f_{i,j}(x,y)=(f_i(x,y), f_j(x,y))$$
\end{defn}
Here $(x,y)$ belongs to $\mathbb{B}_2$.
\noindent 
In the rest of article, we shall use $IVT_{i,j}$ instead of $IVT^{2,2}_{i,j}$ for notational simplicity.\\

\noindent
Here we define a discrete dynamical system over $\mathbb{N} \times \mathbb{N}$ as $$x_{n+1}=IVT_{i,j}(x_n)$$ where $x_n \in \mathbb{N} \times \mathbb{N}$.  

\begin{defn}
	For every $x_0 \in \mathbb{N} \times \mathbb{N}$, if trajectory of the discrete dynamical system $x_{n+1}=IVT_{i,j}(x_n)$ eventually converges to a point $(0,0) \in \mathbb{N} \times \mathbb{N}$, then we call such dynamical system $x_{n+1}=IVT_{i,j}(x_n)$ as \emph{Collatz-like IVT}. 
\end{defn}

\noindent
Let $\mathbb{S}$ be the set of all possible two variable Boolean functions. Clearly, there are exactly $16$ such functions.

\noindent
Let $\mathbb{T}$ be the set of all possible $IVT_{i,j}$ defined over the set $\mathbb{N} \times \mathbb{N}$. Clearly, there are exactly $256$ such IVTs. Our aim is to comprehend the dynamics of these maps. 

\noindent
The following properties and remarks in two dimensional integral value transformation from algebraic perspective are given.

\begin{enumerate}
\item {\bf Properties:} 
\begin{enumerate}[a]
\item The space $(\mathbb{B}_2, +_2, *_2) $ forms a vector space over $\mathbb{F}_2=\{0, 1\}$. Here $ +_2$ and $*_2$ are modulo 2 addition and scalar multiplication. The dimension of the vector space is $2$.

\item The basis of the vector space $(\mathbb{B}_2, +_2, *_2) $ consists of two elements $(0,1)$ and $(1,0)$. The dimension of the vector space is $2$.

\item The space $(\mathbb{S}, \oplus_2, \otimes_2) $ forms a vector space over $\mathbb{F}_2=\{0, 1\}$. Here $\oplus_2$, and $\otimes_2$ are modulo $2$ point wise addition and scalar point-wise multiplication. 

\item The basis of the vector space $(\mathbb{S}, \oplus_2, \otimes_2) $ consists of $4$ elements $f_{1}, f_2, f_{4}, f_{8}$. The dimension of the vector space is $4$.

\item The space $(\mathbb{T}, \oplus_2, \otimes_2) $ forms a vector space over $\mathbb{F}_2=\{0, 1\}$. Here $\oplus_2$, and $\otimes_2$ are modulo 2 function addition and scalar function multiplication.

\item 	The basis of the vector space $(\mathbb{T}, \oplus_2, \otimes_2) $ consists of $16$ elements $IVT_{1,1},$ $IVT_{1,2},$ $IVT_{1,4},$ $IVT_{1,8},$ $IVT_{2,1},$ $IVT_{2,2},$ $IVT_{2,4},$ $IVT_{2,8},$ $IVT_{4,1},$ $IVT_{4,2},$ $IVT_{4,4},$ $IVT_{4,8},$ $IVT_{8,1},$ $IVT_{8,2},$ $IVT_{8,4},$ $IVT_{8,8}$. The dimension of the vector space is $16$.

\item There are $4$ linear Boolean functions viz. $f_0, f_6, f_{10}$ and $f_{12}$ in $\mathbb{S}$. 

\item There are $16$ linear maps in $\mathbb{T}$ and they are $IVT_{0,0},$ $IVT_{0,6},$ $IVT_{0,10},$ $IVT_{0,12},$ $IVT_{6,0},$ $IVT_{6,6},$ $IVT_{6,10},$ $IVT_{6,12},$ $IVT_{10,0},$ $IVT_{10,6},$ $IVT_{10,10},$ $IVT_{10,12},$ $IVT_{12,0},$ $IVT_{12,6},$ $IVT_{12,10}$ and $IVT_{12,12}$. 

\item 	Sum and Composition to two linear maps in $\mathbb{T}$ is linear.

\item There are $24$ bijective $f_{i,j}$ in $\mathbb{S}$. The twelve of the bijective functions are $f_{3,5}, f_{3,6}, f_{3,9},$ $f_{3,10}, f_{5,6}, f_{5,9}, f_{5,12},$ $ f_{6,10}, f_{6,12}, f_{9,10}, f_{9,12},$ and $f_{10,12}$.

\item There are only $6$ isomorphisms (linear and bijective) exist in $\mathbb{T}$. They are $f_{6,10},$ $f_{6,12},$, $f_{10,6},$ $f_{10,12},$ $f_{12,6},$ and $f_{12,10}$.      

\end{enumerate}

\item {\bf Remarks:}
\begin{enumerate}[a]
\item 	$IVT_{i,j} \in \mathbb{T}$ is a basis element of $\mathbb{T}$ if and only if $f_i$ and $f_j$ are also basis elements in $\mathbb{S}.$

\item $IVT_{i,j} \in \mathbb{T}$ is linear if and only if $f_i$ and $f_j$ are also linear.

\item If $f_{i,j}$ is linear then $f_{j,i}$ is also linear.

\item None of the basis function in $\mathbb{T}$ is linear.

\item If $f_{i,j}$ is bijective then $f_{j,i}$ is also bijective.
\end{enumerate}
\end{enumerate}

\subsection{IVTs through STDs}
\noindent
In graph theory, a directed graph $G$ is defined as $G(V,E,F)$; where $F$ is a function such that every edge in $E$ an oder pair of vertices in $V$ \cite{28}. In our context, we have designed the \emph{state transition diagram} (STD) of $f_{i,j}$ as directed graph. 

{\bf Illustration:} For an example, consider the domain of $2$-variable Boolean functions which have four unique possibilities i.e. $00$, $01$, $10$ and $11$ for its constituent variables as shown in Table 1. 

\begin{table}[!ht]
\label{Tab:1}
\tbl{The truth table definition of the 2-vriable functions $f_{6,13}$ and $f_{0,12}$ and their transitions.}%
{\begin{tabular}{|c|c|c|c|c|c|c|}
\hline
Variables & \multicolumn{2}{c|} {$f_{6,13}$} & Transitions & \multicolumn{2}{|c|} {$f_{0,12}$} & Transitions\\
\cline{1-3} \cline{5-6}

$x_{2}~x_{1}$ & $f_6$ & f$_{13}$ & & $f_0$ & f$_{12}$ & \\
\hline
0~~~0 & 0 & 1 & (0,0)$\rightarrow$(0,1) & 0 & 0 & (0,0)$\rightarrow$(0,0)\\
0~~~1 & 1 & 0 & (0,1)$\rightarrow$(1,0) & 0 & 0 & (0,1)$\rightarrow$(0,0)\\
1~~~0 & 1 & 1 & (1,0)$\rightarrow$(1,1) & 0 & 1 & (1,0)$\rightarrow$(0,1)\\
1~~~1 & 0 & 1 & (1,1)$\rightarrow$(0,1) & 0 & 1 & (1,1)$\rightarrow$(0,1)\\
\hline
\end{tabular}}
\end{table}

\noindent 
From the truth table representation of $2$-variable Boolean function $f_{6,13}$, we find that $(0,0) \rightarrow (0,1)$, $(0,1)\rightarrow (1,0)$, $(1,0)\rightarrow (1,1)$ and $(1,1)\rightarrow (0,1)$ are the required transitions. Similarly, the transitions are shown for the function $f_{0,12}$. Based on these four transitions a directed graph called STD ($G=(V,E)$) can be drawn where the vertices or nodes in the graph are four possible bit pairs i.e. $V= \mathbb{B}_2$ and the directed edges $(E)$ are the transitions or mapping from the set $\mathbb{B}_2$ to $\mathbb{B}_2$. So, STD identifies the local transformations for all possible combinations of domain set to the rules, given the truth table of one or more rules of $n$-variable. The STDs of function $f_{6,13}$ and $f_{0,12}$ are shown in Fig. 1 and it is noted that there will be always four nodes each node with an out degree edge/transition, but the in degree edge for a particular node can vary from zero to maximum of four. From Fig. 1 ($f_{6,13}$), the in-degree of node $(0,1)$ is 2, $(1,1)$ is 1, $(1,0)$ is 1, and $(0,0)$ is 0. Total in-degrees are the summand of all in-degrees of four nodes that must be equal to $4$. 

\begin{figure} [!ht]
\centering
\tcbox[sharp corners, boxsep=0mm, boxrule=.3mm, 
colframe=green!30!black, colback=white]{\includegraphics[]{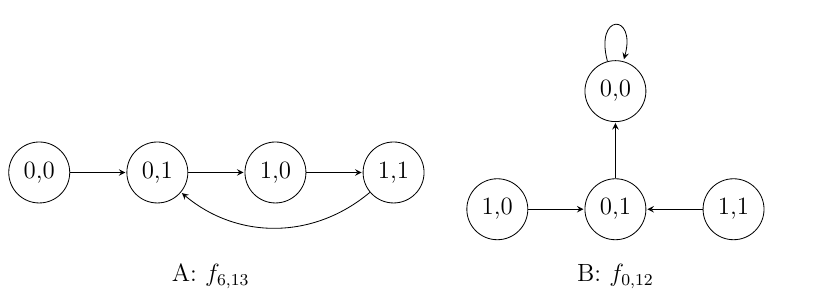} }
\caption{The STDs of the function $f_{6,13}$ (A) and $f_{0,12}$ (B).}
\label{fig:fig1}
\end{figure}

\noindent
The STDs of linear functions are given in Fig. 2. The STDs of only linear $IVT_{i,j}$ are given, we omit other STDs of linear $f_{j,i}$. The STDs of basis functions are given in Fig. 3. The STDs of only basis $IVT_{i,j}$ are given, we omit other STDs of basis $f_{j,i}$. The STDs of bijective functions are given in Fig. 4. The STDs of only bijective $IVT_{i,j}$ are given, we skip STDs of other bijective $f_{j,i}$. The STDs of isomorphisms are given in Fig. 5. The STDs of only $IVT_{i,j}$ which are isomorphic are given, we skip other STDs of isomorphism $f_{j,i}$.

\begin{figure} [!ht]
	\centering
	\tcbox[sharp corners, boxsep=0mm, boxrule=.3mm, 
	colframe=green!30!black, colback=white]{\includegraphics[]{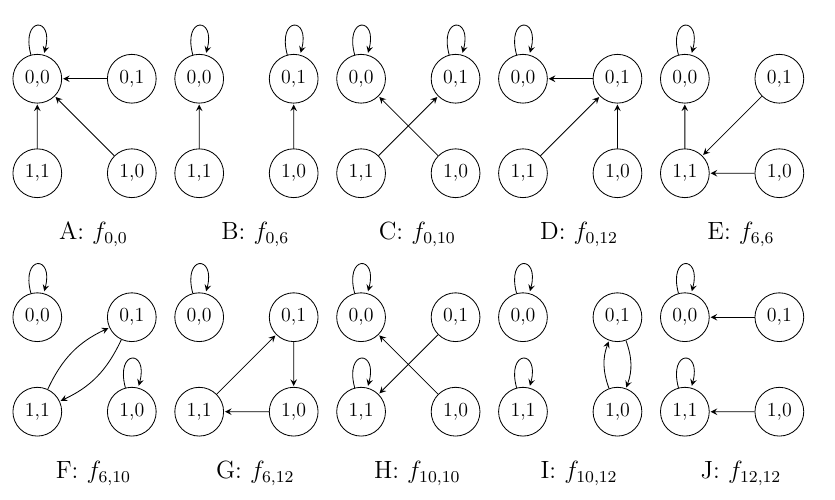} }
	\caption{The STDs of 10 linear functions $f_{i,j}$.}
	\label{fig:fig2}
\end{figure}

\begin{figure} [!ht]
	\centering
	\tcbox[sharp corners, boxsep=0mm, boxrule=.3mm, 
	colframe=green!30!black, colback=white]{\includegraphics[]{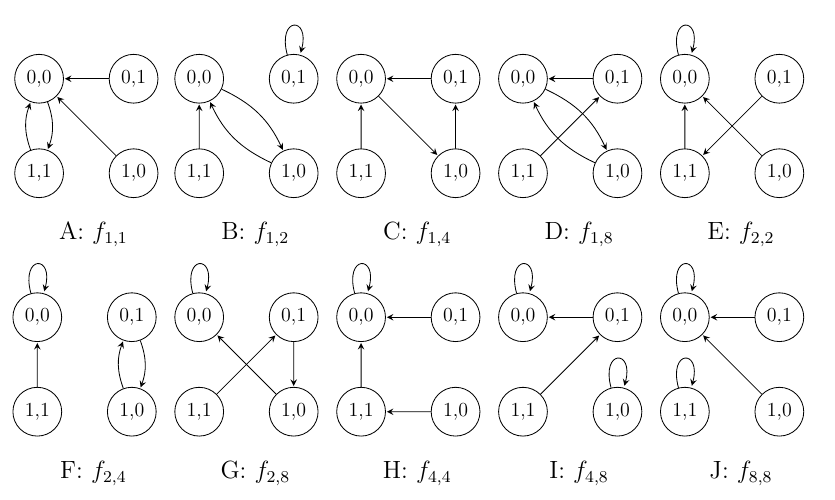}} 
	\caption{The STDs of 10 basis functions $f_{i,j}$.}
	\label{fig:fig3}
\end{figure}

\begin{figure} [!ht]
	\centering
	\tcbox[sharp corners, boxsep=0mm, boxrule=.3mm, 
	colframe=green!30!black, colback=white]{\includegraphics[]{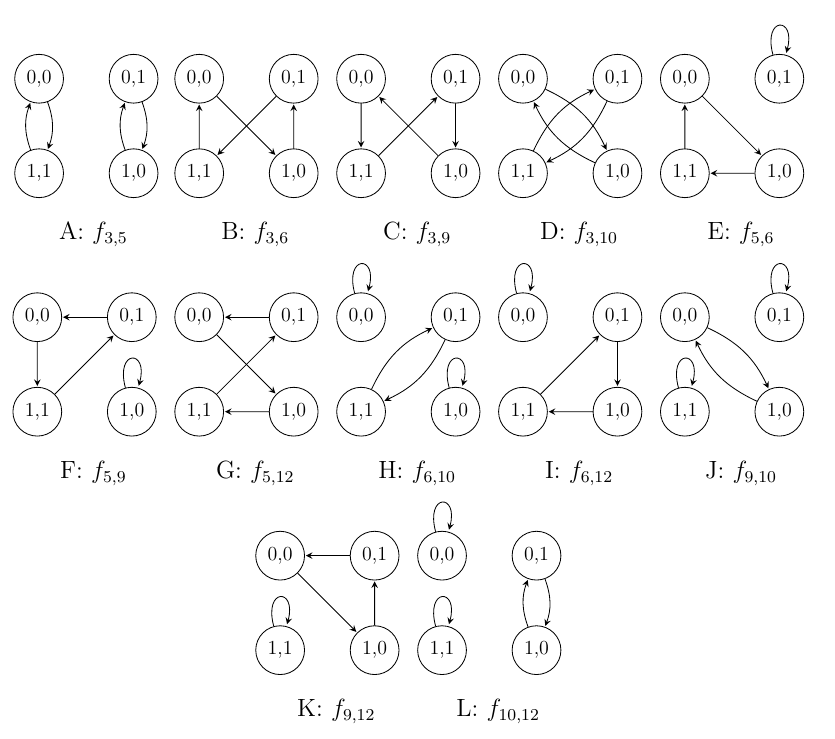} }
	\caption{The STDs of 12 bijective functions $f_{i,j}$.}
	\label{fig:fig4}
\end{figure}

\begin{figure} [!ht]
	\centering
	\tcbox[sharp corners, boxsep=0mm, boxrule=.3mm, 
	colframe=green!30!black, colback=white]{\includegraphics[]{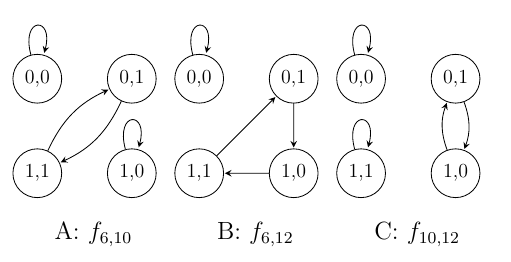}} 
\caption{The STDs of 3 isomorphisms functions $f_{i,j}$.}
	\label{fig:fig5}
\end{figure}

\section{Dynamics of $IVT_{i,j}$ Using State Transition Diagram of $f_{i,j}$}
\subsection{Dynamics of STDs when p=2, k=2}
\noindent
For $i, j \in \{0, 1,\dots,15\}$ in two dimensions, there are $2^{4}\times 2^{4}=256$ functions $f_{i,j}$. From the STDs of those $256$ functions $f_{i,j}$, their dynamics are pretty clear. There are specifically two kinds of dynamics possible; one is attracting to a single attractor and other kind is periodic attractor of period ($pr$), $2 \leq pr \leq 4$. For examples, the single attractor function is $f_{0,12}$ and period $3$ attractor function is $f_{6,13}$ (Fig. 1). 

\begin{pro}
The dynamical system $x_{n+1}=IVT_{i,j}(x_n)$ is \emph{Collatz-like} if and only if STD of $f^k_{i,j}((x,y))=(0,0)$ for all $(x,y) \in \mathbb{B}_2$, for some $k \in \mathbb{N}$. Here $f^k_{i,j}$ denotes the k-composition of the function defined as $f^k_{i,j}=f_{i,j}(f^{k-1}_{i,j})$. 
\end{pro}

\begin{thm}
There are exactly $16$ \emph{Collatz-like} $x_{n+1}=IVT_{i,j}(x_n)$ in $\mathbb{T}$.
\end{thm}

\begin{proof}
As per definition of Collatz-like function, we are looking for functions $f_{i,j}$ such that for some $k$, $$f^k_{i,j}((x,y))=(0,0)$$ for all $(x,y) \in \mathbb{B}_2$. We will check and count through the possible variety of STDs structure of $f_{i,j}$.

Suppose $f_{i,j}$ has the required property and let $G_{f_{i,j}}$ be the directed graph (STD) with four vertices- $a$ is $(0,0)$, $b$ is $(0,1)$, $c$ is $(1,0)$, $d$ is $(1,1)$ and an edge $(x,y) \to (m,n)$ whenever $f_{i,j}(x,y) = (m,n)$.  Each vertex of $G_{f_{i,j}}$ has out degree 1. We make the following claims, all of which are simple.

\noindent
1. $G_{f_{i,j}}$ must have at least one cycle, say $C$.\\
2. $G_{f_{i,j}}$ must be connected, and it has only one component.\\
3. The cycle $C$ is the only cycle in $G_{f_{i,j}}$, since every vertex must be in the same component as $(0,0)$.  \\
4. The cycle $C$ is a sink: for each $(x,y)$, $f^k_{i,j}(x,y)$ is in $C$ for all sufficiently large $k$ and $(0,0) \in C$.\\
5. If $(0,0) \ne f^k_{i,j}((0,0))$, then $f^k_{i,j}(a)\ne f^k_{i,j}(f((0,0)))$ for every $k$. And no other vertex is in $C$.\\
	
\noindent
So we have $(0,0) = f_{i,j}((0,0))$. Therefore, we could have only $4$ topologies (STDs) for $G_{f_{i,j}}$ as shown below. The vertex $a$ is fixed and by interchanging $b$, $c$ and $d$ in different levels/positions we can easily enumerate how many elements that have $f^k_{i,j}((x,y)) = (0,0)$.\\

\begin{center}
\begin{tabular}{p{2.2cm}p{3.2cm}p{3.2cm}p{3.2cm}p{3.1cm}} 
Topology: &  
\begin{tikzpicture}[
> = stealth, % arrow head style
shorten > = .4pt, % don't touch arrow head to node
auto,
node distance = 1cm, % distance between nodes
semithick % line style
]	
\tikzstyle{every state}=[
draw = black,
thick,
fill = white,
minimum size = 1mm
]
\node[state] (1) {a};
\node[state] (2) [below of=1] {c};
\node[state] (3) [left of=2] {b};
\node[state] (4) [right of=2] {d};

\path[->] (1) edge [loop above] node {} (1);
\path[->] (2) edge node {} (1);
\path[->] (3) edge node {} (1);
\path[->] (4) edge node {} (1);
\end{tikzpicture}
&
\begin{tikzpicture}[
> = stealth, % arrow head style
shorten > = .5pt, % don't touch arrow head to node
auto,
node distance = 1cm, % distance between nodes
semithick % line style
]	
\tikzstyle{every state}=[
draw = black,
thick,
fill = white,
minimum size = 2mm
]
\node[state] (1) {a};
\node[state] (2) [below left of=1] {b};
\node[state] (3) [below right of=1] {c};
\node[state] (4) [below of=2] {d};

\path[->] (1) edge [loop above] node {} (1);
\path[->] (2) edge node {} (1);
\path[->] (3) edge node {} (1);
\path[->] (4) edge node {} (2);
\end{tikzpicture} 
& 
\begin{tikzpicture}[
> = stealth, % arrow head style
shorten > = .5pt, % don't touch arrow head to node
auto,
node distance = 1cm, % distance between nodes
semithick % line style
]	
\tikzstyle{every state}=[
draw = black,
thick,
fill = white,
minimum size = 2mm
]
\node[state] (1) {a};
\node[state] (2) [below of=1] {b};
\node[state] (3) [below left of=2] {c};
\node[state] (4) [below right of=2] {d};

\path[->] (1) edge [loop above] node {} (1);
\path[->] (2) edge node {} (1);
\path[->] (3) edge node {} (2);
\path[->] (4) edge node {} (2);
\end{tikzpicture} 
&
\begin{tikzpicture}[
> = stealth, % arrow head style
shorten > = .3pt, % don't touch arrow head to node
auto,
node distance = 1cm, % distance between nodes
semithick % line style
]	
\tikzstyle{every state}=[
draw = black,
thick,
fill = white,
minimum size = 1mm
]
\node[state] (1) {a};
\node[state] (2) [below of=1] {b};
\node[state] (3) [below of=2] {c};
\node[state] (4) [below of=3] {d};

\path[->] (1) edge [loop above] node {} (1);
\path[->] (2) edge node {} (1);
\path[->] (3) edge node {} (2);
\path[->] (4) edge node {} (3);
\end{tikzpicture} \\ \hline
Number of functions: & 1 & 6 & 3 & 6 \\ \hline
Functions ($f_{i,j}:(i,j)$): & (0,0) &  (0,4), (0,8), (2,0), (2,2), (4,4), (8,0) & (0,12), (6,6), (10,0) & (2,6), (2,8), (4,12), (6,4), (8,4), (10,2) \\ \hline
\end{tabular}
\end{center}
\noindent
Hence, the total number of Collatz-like functions in two dimensions is 16.
\end{proof}

\noindent
It is noted that the attractor $(0,0)$ is a global attractor for the \textit{Collatz-like} IVTs, that is for every initial value taken from $\mathbb{N} \times \mathbb{N}$, the dynamical system takes all the trajectories to $(0,0)$.

\begin{pro}
	The dynamical system $x_{n+1}=IVT_{i,j}(x_n)$ possesses multiple attractor(s) of the form $(2^s-1,0)$, $s\in \mathbb{N}$ if and only if $f^k_{i,j}((x,y))=(1,0)$, for all $(x,y) \in \mathbb{B}_2$, for some $k \in \mathbb{N}$. 
\end{pro}

\begin{thm}
	There are countably infinite number of local attractors of the form $(2^s-1,0)$, $s\in \mathbb{N}$. 
\end{thm}

\begin{proof}
	For any $x_0 \in \mathbb{N} \times \mathbb{N}$, the dynamical system eventually proceeds to a point of the form $(2^s-1,0)$, according to definition of the dynamical system map $IVT_{i,j}$. It is obvious that the set of all possible points of the form $(2^s-1,0)$, is a subset of the $\mathbb{N} \times \mathbb{N}$. Since the set $\mathbb{N} \times \mathbb{N}$ is countable and so the subset is countable. Hence the theorem is followed.   
\end{proof}

\begin{pro}
	The dynamical system $x_{n+1}=IVT_{i,j}(x_n)$ possesses multiple attractor of the form $(0, 2^s-1)$ if and only if $f^k_{i,j}((x,y))=(0,1)$, for all $(x,y) \in \mathbb{B}_2$, $s\in \mathbb{N}$ for some $k \in \mathbb{N}$.
\end{pro} 

\begin{pro}
	There are countably infinite number of local attractors of the form $(0,2^s-1)$, $s\in \mathbb{N}$ for the dynamical system $x_{n+1}=IVT_{i,j}(x_n)$. 
\end{pro}

\begin{pro}
	The dynamical system $x_{n+1}=IVT_{i,j}(x_n)$ possesses multiple attractors of the form $(2^s-1,2^t-1)$, $s, t\in \mathbb{N}$ if and only if $f^k_{i,j}=(1,1)$ for some $k \in \mathbb{N}$.
\end{pro} 

\begin{pro}
	There are countably infinite number of local attractors of the form $(2^s-1,2^t-1)$ $s, t\in \mathbb{N}$ for the dynamical system $x_{n+1}=IVT_{i,j}(x_n)$. 
\end{pro}

\noindent
It is noted that there are STDs of some function $f_{i,j}$ such that there are periodic attractor of period $p \geq 2$. \\ In such cases, the dynamical system $x_{n+1}=IVT_{i,j}(x_n)$ will possess to multiple attractors of form $(u,2^t-1), (2^t-1,u), (0,u), (u,0)$ and  $(u,v)$. There are few dynamical systems which are sensitive to the initial conditions. In the following section, we shall classify all such types of dynamical systems based on their dynamical behavior.

\subsection{STD dynamics when $p\geq 2$, $k\geq3$}
In $p$-adic and $k$-variable functions (or $k$-dimension), we have $p^{p^{k}}$ number of functions $f_{i}$ where $i\in \{0,1,\cdots,2^{k}-1\}$. Depending on $k$ value, we select the number of functions in $k$ dimensions and in combining we represent as $f_{i_{1},i_{2},\cdots ,i_{k}}$ and length of each function is $p^{k}$.

In general, for any $i, j \in \{0, 1,\dots2^{k}-1\}$ in $k$ dimensions, there are $p^{p^{k}}\times p^{p^{k}}=p^{2p^{k}}$ functions $f_{i_{1},i_{2},\cdots ,i_{k}}$. In that case, we could have the attractor of period ($pr$), $1 \leq pr \leq p^{k}$. The number of functions is exponentially increasing with the increment of $p$ and $k$. In principle, the discrete dynamics of such IVTs can also be obtained using STDs as we did it in two dimensions $(p=2,k=2)$. The vertex of the STD in $p$ and $k$ dimension is $(a_{1},a_{2},\cdots a_{k})$; where $a_{i}\in \{0,1,\cdots p^{k}\}$. For an example, the STDs of two \textit{Collatz-like} functions $f_{128,72,160}$ ($p=2,k=3$) and $f_{2187,12690}$ ($p=3,k=2$) and a STD of function $f_{41,19230}$ ($p=3,k=2$) with attractors 5 and 2 is shown in Fig. 6.

\begin{figure} [!ht]
\centering
\tcbox[sharp corners, boxsep=0mm, boxrule=.31mm, 
colframe=green!30!black, colback=white]{\includegraphics[]{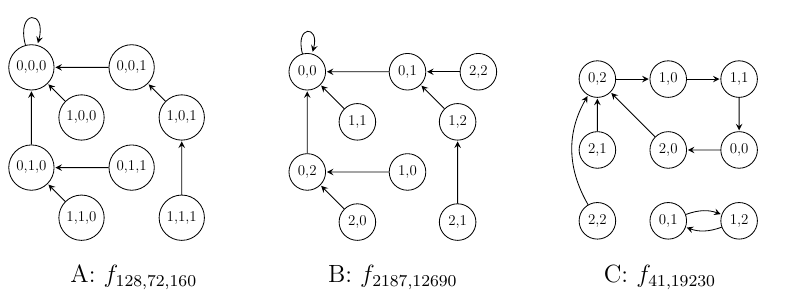} }
\caption{The STD of $f_{128,72,160}$ ($p=2,k=3$), $f_{2187,12690}$ and $f_{41,19230}$ ($p=3,k=2$).}
\label{fig:fig6}
\end{figure}

\noindent
The number of global attractors for different $p$ and $k$ value is shown in Table 2.

\begin{table}[!ht]
\tbl{Number of global attractor for different $p$ and  $k$ values.}%
{\begin{tabular}{|c|c|c|}
\hline
Base ($p$) & Dimension ($k$) & Number of Attractors \\
\hline
p=2 & k=2 & 16 \\ \hline
p=2 & k=3 & 262144 \\ \hline
p=3 & k=2 & 4782969 \\ \hline
\end{tabular}}
\end{table}

Based on the computational evidence, the total number of \textit{Collatz-like} IVTs in the base $p$ with $k$ dimension is conjectured.

\begin{conj}
The total number of Collatz-like IVTs in the base $p$ with $k$ dimension is $p^{k(p^{k}-2)}$.
\end{conj}

\begin{pro}
If $l_{k}$ defines the path length from a particular point to the global attractor in the dynamics of $IVT^{p,k}_{i,j}$, then $1\leq l_{k}\leq (p^{k}-1)$.
\end{pro}

\section{Dynamics of $x_{n+1}=IVT_{i,j}(x_n)$}
Before getting into the detailed classifications of dynamics, a vivid example of a dynamical system $x_{n+1}=IVT_{i,j}(x_n)$, $x_n \in \mathbb{N} \times \mathbb{N}$ is adumbrated with different initial positive natural integers and their trajectories as shown in the Fig. 7. For an example, $IVT_{13,3}(0,2)=(1,3)$ i.e. $(0,2)\rightarrow(1,3)$. The trajectories (path) can be changed with the changes of $p$ value even if the taken function is fixed.

\begin{figure}[!ht]
\begin{center}
\tcbox[sharp corners, boxsep=0mm, boxrule=.35mm, 
colframe=green!30!black, colback=white]{\includegraphics[scale=1]{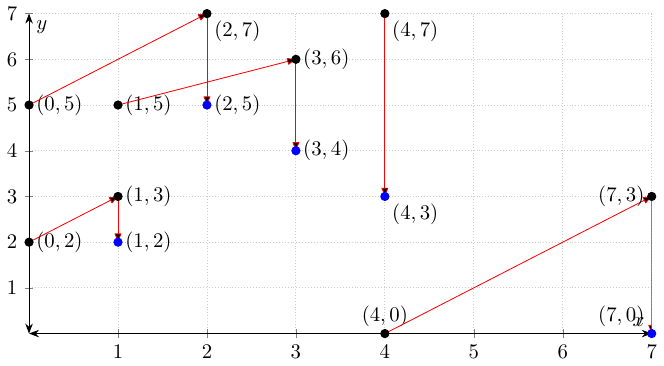}}
\end{center}
\caption{The dynamical system of $IVT_{13,3}$ for the set of initial integer pairs.}
\label{fig:fig7}
\end{figure}

\subsection{Classifications of $IVT_{i,j}$ based on their dynamics when $p=2$, $k=2$}
In this section, we have classified the $256$ of $IVT_{i,j}$ into four disjoint classes ($Class~I-IV$) according to the attractors of four different cycle length(s). For every class, several subclasses are identified according to the dynamic behavior of $IVT_{i,j}$ over $\mathbb{N}\times \mathbb{N}$. It can be seen that there are many $(x,y)$ components in the dynamics of  $IVT_{i,j}$  where $(x,y)\in \mathbb{N}\times \mathbb{N}$. For the simplicity of graph visualization, we have taken certain limit for both $x$ and $y$ i.e. $x,y\in\{0,1\cdots 7\}$. For some of the $IVTs$, all the components are connected (e.g. all $Collatz-like$ IVTs) and other cases there could be of different variety.  For understanding the complete dynamics over the set of positive natural numbers in base $p$ and dimension $k$, $x,y\in \{0,1,\cdots, p^{k}-1\}$ that covering all possible transitions in their $p$-adic representations. Here, $p^{k}-1$ is the minimum value that can be extended to desire limit.

\subsubsection{Class-I: Dynamics of $IVT_{i,j}$: attractor with a fixed point}
A total of 125 $IVT_{i,j}$ are classified into the \textit{Class-I} as shown in Table 3. It is noted that each of the IVTs in the \textit{Class-I} is having attractor with one length cycle (singleton attractor). There are 16 IVTs which are Collatz-like and rest are non-Collatz-like. For the non-Collatz-like dynamical systems, there are countable many local attractors with one cycle.

\begin{table}[!ht]
\centering
\tbl{Dynamics of $IVT_{i,j}$ with singleton attractor or a fixed point.} 
{\begin{tabular}{|p{1.1cm}|p{9cm}|p{2.2cm}|p{1.8cm}|}
\hline
Subclass & $IVT_{i,j}: (i,j)$ & Attractor(s) & Cardinality   \\
\hline
1 & (0,0), (0,4), (0,8), (0,12), (2,0), (2,2), (2,6), (2,8), (4,4), (4,12), (6,4), (6,6), (8,0), (8,4), (10,0), (10,2)
&  $(0,0)$ & 16 (Collatz-like)  \\ \hline

2 & (0,3),  (0,7),  (0,11), (0,15), (1,14),  (1,6), (1,11), (1,15), (4,7),  (4,15), (5,14), (5,15), (8,3), (8,7), (9,6), (9,7) &  $(0,2^n-1)$ & 16 \\ \hline

3 & (5,0), (5,8), (6,1), (6,9), (7,0), (7,2), (7,8), (7,9), (13,0), (13,1), (14,1), (14,3), (15,0), (15,1), (15,2), (15,3) & $(2^n-1,0)$ & 16 \\ \hline

4 & (9,9), (9,13), (10,11), (10,15), (11,9), (11,11), (11,14), (11,15), (13,12), (13,13), (14,13), (14,15), (15,12), (15,13), (15,14), (15,15) & $(2^n-1,2^n-1)$ & 16 \\  \hline

5 & (8,11), (8,15), (9,11), (9,14), (9,15), (12,15), (13,14), (13,15)  & $(x,2^n-1)$ & 8 \\
\hline

6 & (13,8), (13,9), (14,9), (14,11), (15,8), (15,9), (15,10), (15,11)   & $(2^n-1,x)$ & 8 \\ \hline

7 & (0,2), (0,6), (0,10), (0,14), (4,6), (4,14), (8,2), (8,6) & $(0,x)$ & 8 \\ \hline 
8 & (4,0), (4,8), (6,0), (6,2), (6,8), (12,0), (14,0), (14,2) & $(x,0)$ & 8 \\ \hline 

9 & (8,8), (8,12), (10,8), (10,10), (10,14), (12,12),(14,12), (14,14)  & $(x,x)$ & 8 \\ \hline 

10 & (4,11), (5,10), (5,11), (12,3), (13,2)  & $(x,y)$ & 5 \\ \hline 

11 &  (4,2), (4,3), (4,10),  (5,2), (8,10), (8,14), (12,2), (12,8), (12,10), (12,11), (12,14), (13,3), (13,10), (13,11), (14,8), (14,10) & Sensitive to the initial condition & 16 \\ \hline 
\end{tabular}}
%	}
\end{table}

A dynamical system  $IVT_{10,0}$ having an attractor with one length cycle is given in the Fig. 8.

\begin{figure}[!ht]
	\begin{center}
	\tcbox[sharp corners, boxsep=0mm, boxrule=.3mm, 
	colframe=green!30!black, colback=white]{\includegraphics[scale=0.5]{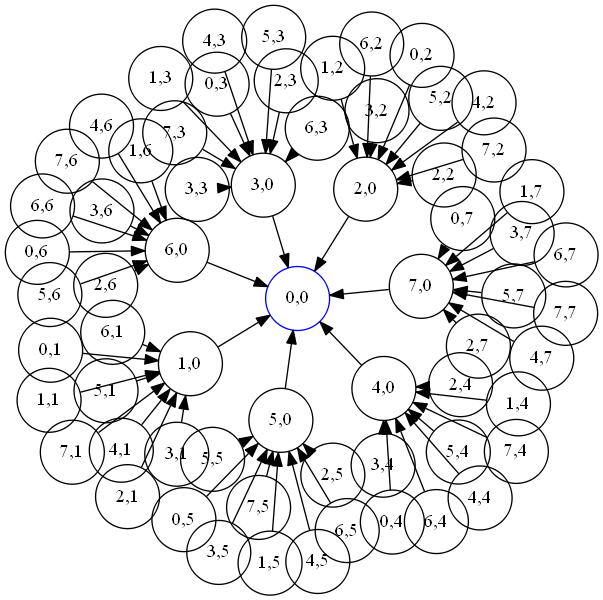}}
	\end{center}
	\caption{The dynamical system of $IVT_{10,0}$ having an attractor with  length  one cycle. Here a set of 64 initial points from $\mathbb{N}\times \mathbb{N}$ are taken, which are in-circled. The point $(0,0)$ is the global attractor.}
	\label{fig:fig8}
\end{figure}

\noindent
There are other 16 out of 125 which are sensitive to the initial conditions that is their trajectories are depending on initial values. Six possible $IVT_{i,j}$ are taken and for all the cases five integer pairs trajectories are shown in Fig. 9. It is noted that each of these initial condition (taken pair) provides sensitive dynamical system ending with different cyclic attractors of length one.

\begin{figure}[!ht]
\begin{center}
\includegraphics[scale=.78]{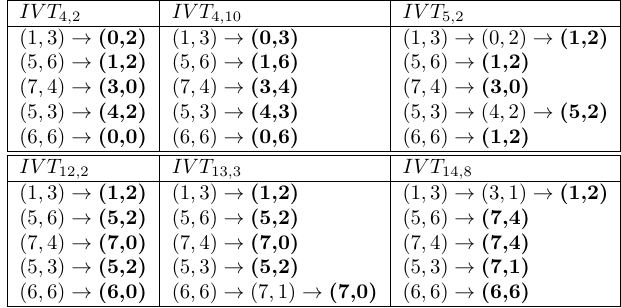}
\end{center}
\caption{Trajectories of some initial values for couple of examples of $IVT_{i,j}$ which are sensitive to the initial conditions. Attractors are marked in bold.}
\label{fig:fig9}
\end{figure}

\subsubsection{Class-II: Dynamics of $IVT_{i,j}$: attractor with two length cycles}

A total of 93 $IVT_{i,j}$ are classified into the \textit{Class-II} as shown in Table 4. All these IVTs in the \textit{Class-II} are having attractors with two length cycles. Some of these dynamical systems are having global attractor of two length cycles. There are 42 out of 93 which are sensitive to the initial conditions that is their trajectories are depending on initial values. Thus, each of these initial conditions provides sensitive dynamical system ending with two length cycle attractor. 

\begin{table}[!ht]
\centering
\tbl{Dynamics of $IVT_{i,j}$ with attractor with two length cycles.} 
{\begin{tabular}{|p{1.1cm}|p{8.5cm}|p{2.5cm}|p{2cm}|}
\hline
Subclass & $IVT_{i,j}: (i,j)$ & Attractor(s) & Cardinality  \\
\hline
1 & (0,1), (0,5), (0,9), (0,13), (4,5), (4,13), (8,1), (8,5)
&  $(0,0)\Leftrightarrow (0,1)$ & 8  \\
\hline
2 & (1,0), (1,8), (3,0), (3,2), (3,8), (9,0), (11,0), (11,2) 
&  $(0,0)\Leftrightarrow (1,0)$ & 8  \\
\hline
3 & (1,1), (1,5), (3,1), (3,3), (3,7), (5,5), (7,5), (7,7)
&  $(0,0)\Leftrightarrow (1,1)$ & 8  \\
\hline
4 & (2,5), (2,13), (3,4), (3,12), (3,13), (10,5), (11,4), (11,5) &  $(x,\bar{x})\Leftrightarrow (\bar{x},x)$ & 8  \\
\hline
5 & (2,11), (2,15), (3,11), (3,14), (3,15), (6,15), (7,14), (7,15) &  $(x,2^{n}-1)\Leftrightarrow (\bar{x},2^{n}-1)$ & 8  \\
\hline
6 & (13,4), (13,5), (14,5), (14,7), (15,4), (15,5), (15,6), (15,7) &  $(2^{n}-1,x)\Leftrightarrow (2^{n}-1,\bar{x})$ & 8  \\
\hline
7 &  (3,5), (3,10), (12,5) &$(p,q)\Leftrightarrow (r,s)$ & 3 \\  \hline
8 &  (1,2), (1,3), (1,7), (1,10), (2,4), (2,10), (2,12), (2,14), (4,1), (4,9), (5,1), (5,3), (5,7), (6,10), (6,11), (6,14), (7,1), (7,3), (7,10), (7,11), (8,9), (8,13), (9,2), (9,8), (9,10), (10,4), (10,12), (10,13), (11,8), (11,10), (11,12), (11,13), (12,1), (12,4), (12,6), (12,7), (12,9), (12,13), (13,6), (13,7), (14,4), (14,6) & Sensitive to the initial condition & 42 \\  \hline
\end{tabular}}
%	}
\end{table}

A dynamical system  $IVT_{4,5}$ having an attractor with two length cycle (global attractor) is given in the Fig. 10. 

\begin{figure}[!ht]
	\begin{center}
	\tcbox[sharp corners, boxsep=0mm, boxrule=.3mm, 
	colframe=green!30!black, colback=white]{\includegraphics[scale=.5]{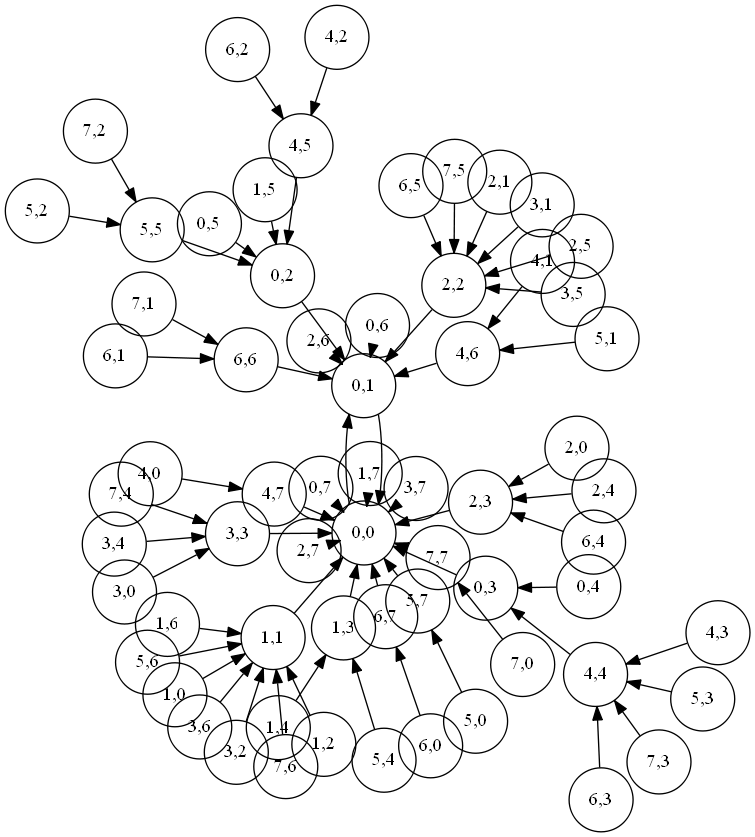}}
	\end{center}
	\caption{The dynamical system of $IVT_{4,5}$ having an attractor with two length cycle.}
	\label{fig:fig10}
\end{figure}

\subsubsection{Class-III: Dynamics of $IVT_{i,j}$: attractor with three length cycles}

A total of 32 $IVT_{i,j}$ are classified into the \textit{Class-III} as shown in Table 5. All these IVTs in the \textit{Class-III} are having attractor with three length cycles. There are some dynamical systems with global three length cycle attractors. There are 8 out of 32 which are sensitive to the initial conditions that is their trajectories are depending on initial values. Each of these initial conditions provides sensitive dynamical systems ending three length cycle attractor.

\begin{table}[!ht]
\centering
\tbl{Dynamics of $IVT_{i,j}$ with attractor with three length cycle.} 
{\begin{tabular}{|p{1.1cm}|p{6.5cm}|p{3cm}|p{2cm}|}
\hline
Subclass & $IVT_{i,j}: (i,j)$ & Attractor(s) & Cardinality  \\
\hline
1 & (1,4),(1,12),(2,1),(2,9),(9,4),(10,1)
&  $(0,0)-(1,0)-(0,1)$ & 6  \\
\hline
2 & (1,9),(1,13),(2,3),(2,7),(5,13),(6,7) 
& $(0,1)-(0,0)-(1,1)$ & 6  \\
\hline
3 & (5,4),(7,4),(7,6),(9,1),(11,1),(11,3)
& $(1,0)-(1,1)-(0,0)$ & 6  \\
\hline
4 & (6,13), (7,12), (7,13), (10,7), (11,6), (11,7)
& $(x,y)-(y,z)-(z,x)$ & 6  \\
\hline
5 & (5,6), (5,9), (6,3), (6,12), (9,3),  (9,12), (10,6), (10,9)
& Sensitive to the initial condition & 8  \\
\hline
\end{tabular}}
\end{table}

A dynamical system  $IVT_{9,1}$ having an attractor with three length cycle is shown in the Fig. 11.

\begin{figure}[!ht]
	\begin{center}
	\tcbox[sharp corners, boxsep=0mm, boxrule=.3mm, 
	colframe=green!30!black, colback=white]{\includegraphics[scale=.45]{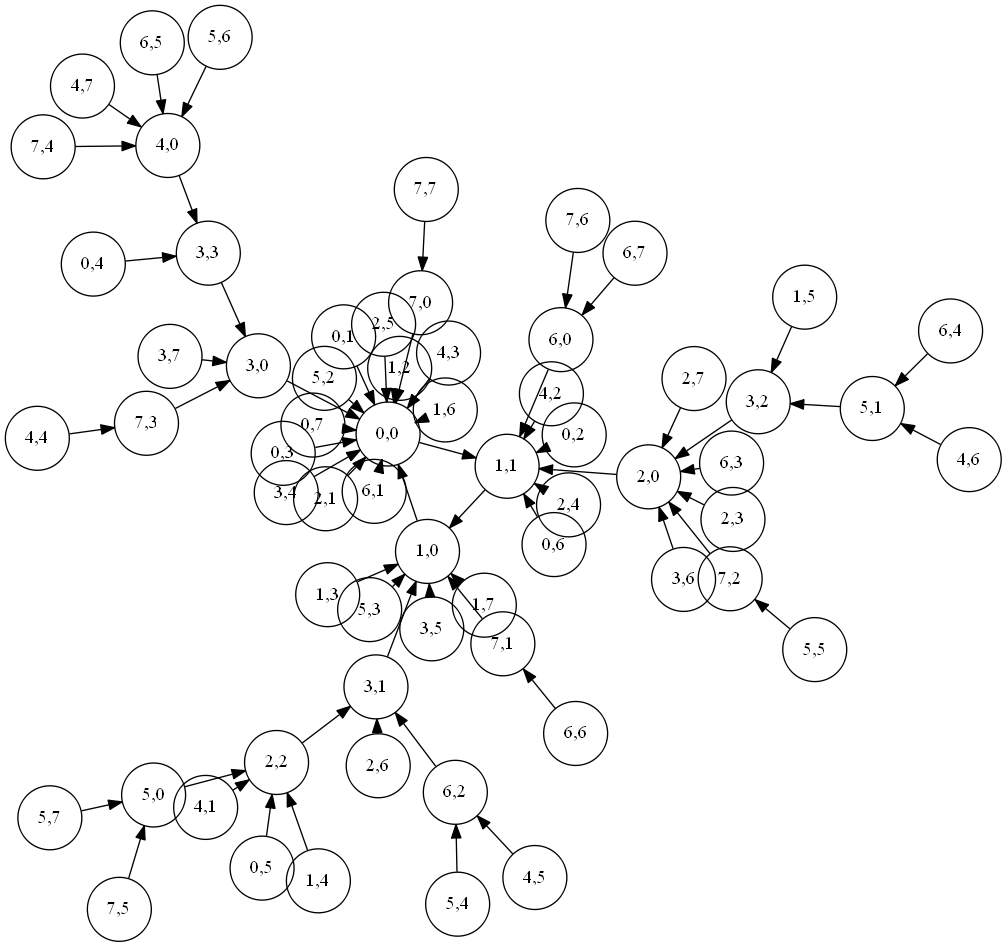}}
	\end{center}
	\caption{The dynamical system  $IVT_{9,1}$ having an attractor with three length cycle.}
	\label{fig:fig11}
\end{figure}

\subsubsection{Class-IV: Dynamics of $IVT_{i,j}$: attractor with four length cycles}

A total of 6 $IVT_{i,j}$ are classified into the \textit{Class-IV} as shown in Table 6. All these IVTs in the \textit{Class-IV} are having attractors with four length cycles. These dynamical systems do possess global attractors with four length cycles. 

\begin{table}[!ht]
\centering
\tbl{Dynamics of $IVT_{i,j}$ with attractor with four length cycles.} 
%\resizebox{12cm}{!}{
{\begin{tabular}{|p{1.1cm}|p{5cm}|p{4cm}|p{2cm}|}
\hline
Subclass & $IVT_{i,j}: (i,j)$ & Attractor(s) & Cardinality \\
\hline
1 & (3,6),(3,9),(5,12),(6,5),(9,5),(10,3)
&  $(0,0)-(0,1)-(1,0)-(1,1)$ & 6  \\
\hline
\end{tabular}}
%	}
\end{table}

A dynamical system  $IVT_{9,5}$ having an attractor with four length cycle is shown in the Fig. 12.

\begin{figure}[!ht]
\begin{center}
\tcbox[sharp corners, boxsep=0mm, boxrule=.3mm, 
colframe=green!30!black, colback=white]{\includegraphics[scale=0.32]{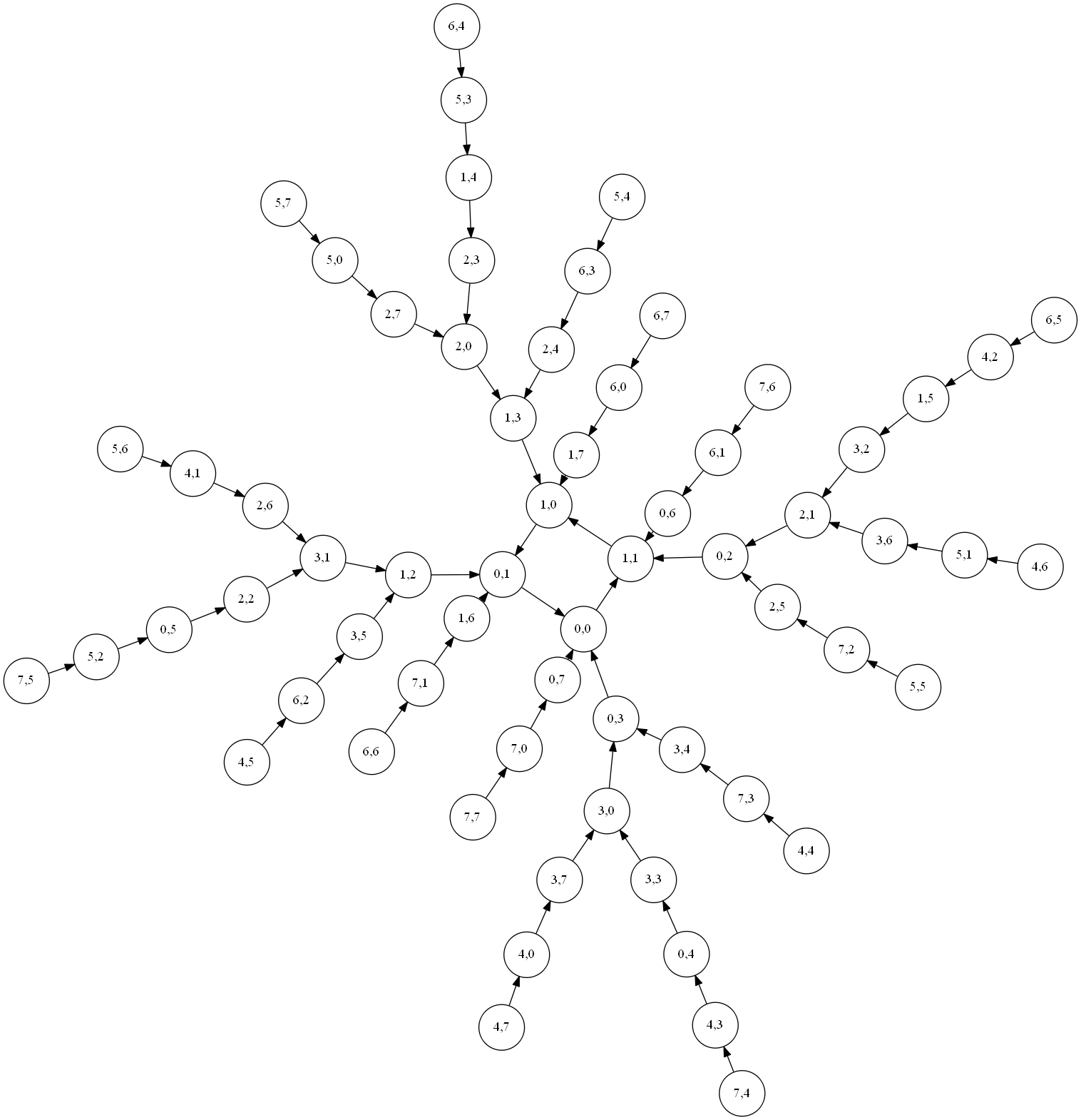}}
\end{center}
\caption{The dynamical system  of $IVT_{9,5}$ having an attractor with four length cycle.}
\label{fig:fig12}
\end{figure}

\subsection{Classifications of $IVT^{p,k}_{i,j}$ based on their dynamics when $p\geq2$, $k\geq 3$}
In case of $p=2$, $k=2$, the cycle length of an attractor observed in the dynamics of $IVT_{i,j}$ does not exceed 4 which is the maximum cycle length of an STD of $f_{i,j}$. It is quite impossible to observe the dynamics of individual $IVT^{p,k}_{i,j}$ when $p$ and $k$ are sufficiently large. From the dynamics of STDs of some $f^{p,k}_{i,j}$, one can classify the corresponding $IVT^{p,k}_{i,j}$ based on their dynamics ends up with cycle length exclusively $1,2,\dots, p^k$ (Fig. 6) into $p^k$ number of classes. The cardinality of different classes based on dynamics of $IVT^{p,k}_{i,j}$ having attractors of different cycle lengths (upto $p^{k}$) including an exclusive outer class are presented for different $p$ and $k$ in Table 7. The exclusive outer class contains $IVT^{p,k}_{i,j}$ such that the length of the cycle of the attractors exceed the maximum cycle length ($p^k$) of the STD of $f^{p,k}_{i,j}$. Such a discrete dynamical systems is called here as \textit{chaotic dynamical system}. \\ 
\noindent
An example of chaotic dynamical system of $IVT^{3,2}_{41,19230}$ using the function $f_{41,19230}$ ($p=3, k=2$) is shown in Fig. 13, whose STD is given in Fig. 6 (C). From the STD of $f_{41,19230}$, we could see two types of cycles (2 and 5) is present. But, the dynamical system shows an additional cycle of length 10 (l.c.m of 2 and 5). The reason can be easily seen from the ternary representation (here, $p=3$) of any pair $(x,y)\in \mathbb{N}\times \mathbb{N}$. Take a pair (5,6) from the cycle of length $10$ from Fig. 13. The ternary representation of $(5,6)$ is $(12,20)$. The pair wise component are $12$ (first component from left) and $20$ (second component from left). These two components are from two different cycles 2 and 5 respectively from STD of $f_{41,19230}$ (Fig. 6(C)). All the STD components are distinct. So, while calculating the IVT one by one iteration, the components $12$ and $20$ will come back/repeat to the same initial stage after  $l.c.m (5,2)$ iterations or cycle of length $l.c.m (5,2)$. This is true for arbitrary pair in any base $p$ and also true if we increase the dimension $k$.

\begin{table}[!ht]
\centering
\tbl{Number of $IVT_{i,j}$ having attractors of different length cycles count in the respective classes.} 
%\resizebox{8cm}{
{\begin{tabular}{|p{1.15cm}|p{1.3cm}|p{1.3cm}|p{1.3cm}|c|p{1.2cm}|c|c|p{1.5cm}|p{1.3cm}|p{1cm}|}
\hline
Category & Class-I & Class-II & Class-III & Class-IV & Class-V & Class-VI & Class-VII & Class-VIII & Class-IX & Outer     \\
\hline
 p=2,k=2 & 125 & 93 &  24 & 6 & $\times$ & $\times$ & $\times$ & $\times$ & $\times$ & $\times$ \\ \hline
 p=2,k=3 & 4782969 & 5752131 & 3009888 & 1692180 & 653184 & 773920 & 46080 &	5040 & $\times$ &	61824 \\ \hline
 p=3,k=2 & 100000000 & 127790505 & 69554240 & 43296120 & 18144000 & 23166360 & 2073600 & 408240 & 40320 & 2947104 \\ \hline
 \end{tabular}}
\end{table}

\begin{figure}[!ht]
\begin{center}
\tcbox[sharp corners, boxsep=0mm, boxrule=.3mm, 
colframe=green!30!black, colback=white]{\includegraphics[scale=0.31]{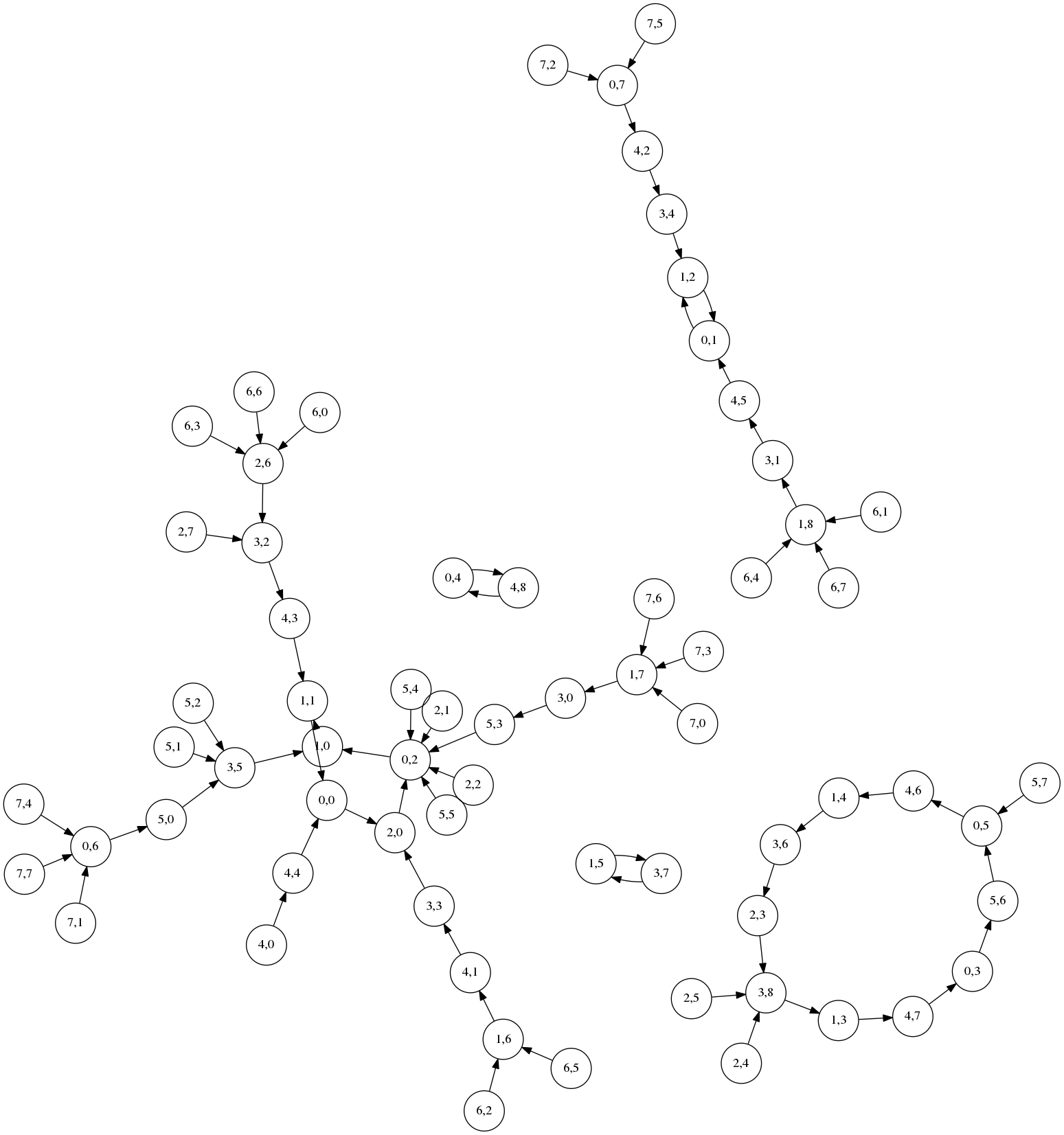}}
\end{center}
\caption{The dynamical system  of $IVT_{41,19230}$ having an attractor with cycle length 2, 5 and 10.}
\label{fig:fig13}
\end{figure}

\begin{pro}
Number of bijective functions ($f^{p,k}_{i_{1},i_{2},\cdots i_{k}}$) in $p$-adic, $k$ dimension is $(p^{k}-1)!.$
\end{pro}

\subsection{Dynamics of two dimensional $IVT_{i,j}$ from algebraic perspective}

Here we go with classification of two dimensional $IVT_{i,j}$ from their algebraic perspectives such as basis, linearity, bijective and isomorphism which are shown in Table 8. 

\begin{table}[!ht]
\centering
\tbl{Dynamics of $IVT_{i,j}$ from Algebraic Perspective} 
%\resizebox{13cm}{!}{
{\begin{tabular}{|p{2.5cm}|p{2cm}|p{4cm}|p{6cm}|}
\hline
Algebraic Character & $IVT_{i,j}: (i,j)$ & Attractor(s) & Remark  \\
\hline
Basis & (1,1) &   $(0,0)\Leftrightarrow (1,1)$ & Period two attractor  \\ \hline

Basis  & (1,2) &   Unstable & Sensitive to the initial condition  \\ \hline

Basis  & (1,4) &   $(0,0)-(1,0)-(0,1)$ & Period three attractor   \\ \hline

Basis  & (1,8) &   $(0,0)\Leftrightarrow (1,0)$ & Period two attractor   \\ \hline

Basis  & (2,2) &   $(0,0)$ & Collatz-like   \\ \hline

Basis  & (2,4) &   Unstable & Sensitive to the initial condition   \\ \hline

Basis  & (2,8) &   $(0,0)$ & Collatz-like   \\ \hline

Basis & (4,4) &   $(0,0)$ & Collatz-like  \\ \hline

Basis & (4,8) &   $(x,0)$ & Period one attractor (Non-Collatz-like)   \\ \hline

Basis  & (8,8) &   $(x,x)$ & Period one attractor (Non-Collatz-like)   \\ \hline

Linear & (0,0) &   $(0,0)$ & Collatz-like   \\ \hline

Linear & (0,6) &   $(0,x)$ & Period one attractor (Non-Collatz-like)    \\ \hline

Linear & (0,10) &   $(0,x)$ & Period one attractor (Non-Collatz-like)    \\ \hline

Linear & (0,12) &   $(0,0)$ & Collatz-like    \\ \hline

Isomorphism & (6,10) &   Unstable & Sensitive to the initial condition    \\ \hline

Isomorphism & (6,12) &   Unstable & Sensitive to the initial condition    \\ \hline

Linear & (6,6) &    $(0,0)$ & Collatz-like    \\ \hline

Isomorphism & (10,12) &   Unstable & Sensitive to the initial condition    \\ \hline

Linear & (10,10) &   $(x,x)$ & Period one attractor (Non-Collatz-like)     \\ \hline

Linear & (12,12) &   $(x,x)$ & Period one attractor (Non-Collatz-like)    \\ \hline

Bijective & (3,5) &  $\lvert p+r \rvert=\lvert q+s \rvert=2^{n}-1$ & Period two attractor \\ \hline

Bijective & (3,6) &  $(0,0)-(0,1)-(1,0)-(1,1)$ & Period four attractor \\ \hline

Bijective & (3,9) &  $(0,0)-(0,1)-(1,0)-(1,1)$ & Period four attractor \\ \hline

Bijective & (3,10) &  $\lvert p+r \rvert=2^{n}-1$ & Period two attractor \\ \hline

Bijective & (5,6) &  Unstable & Sensitive to the initial condition \\ \hline

Bijective & (5,9) &  Unstable & Sensitive to the initial condition  \\ \hline

Bijective & (5,12) &  Unstable & Sensitive to the initial condition  \\ \hline

Bijective & (9,10) &  Unstable & Sensitive to the initial condition   \\ \hline

Bijective & (9,12) &  Unstable & Sensitive to the initial condition   \\ \hline
\end{tabular}}
\end{table}

\noindent
It is to be noted that the dynamics of the $IVT_{i,j}$ is not always same as the dynamics of $IVT_{j,i}$ although they belong to same algebraic class as previously pointed out.

\section{Conclusions and Future Endeavors}
In this manuscript, notion of dimension preservative \textit{Integral Value Transformations} (IVTs) is defined over $\mathbb{N}^k$ using the set of $p$-adic functions. The main focus is given in order to understand the two dimensional integral value transformation ($p=2$ and $k=2$, which are computationally feasible) and studied the dynamics of all possible pair functions in four different kinds of attractors. Among them, a special class of IVTs (Collatz-like) is computationally identified along with some others global and local dynamical systems. In addition, from the algebraic properties of $IVT_{i,j}$, the dynamics are critically observed. From the current analysis, it is evident that as $p$ and $k$ go high, the total numbers of Collatz like IVTs are increasing very fast. Consequently, the dynamical behavior of these IVTs is computationally difficult to simulate and apprehend. But, some of their dynamics can be vivid from our current analysis in two dimensions and some other cases. In this concern, we have shown some quantitative informations for other two special cases of IVTs ($p=2$,$k=3$ and $p=3$, $k=2$). Further, the analysis also clearly evident that as $p$ and $k$ go high, we come across with some chaotic dynamical behavior. In future endeavor, one can explore with some higher $p$ and $k$ considering the limitation of simulation.\\

\subsection*{Acknowledgment}
We thank to anonymous reviewers for their constructive suggestions and valuable directions that help to improve the article quality.

\bibliographystyle{ws-ijbc}
\bibliography{sample}
\end{document}